\documentclass{amsart}
\usepackage{amsfonts}

\usepackage{graphicx}
\usepackage{amscd}

\newtheorem{theorem}{Theorem}
\theoremstyle{plain}

\newtheorem{corollary}{Corollary}

\newtheorem{definition}{Definition}
\newtheorem{example}{Example}

\newtheorem{lemma}{Lemma}

\newtheorem{remark}{Remark}

\numberwithin{equation}{section}

\input{tcilatex}

\begin{document}
\title[On the Cayley digraphs that are patterns of unitary matrices]{On the Cayley digraphs that are \\
patterns of unitary matrices}
\author{Simone Severini}
\address{University of Bristol, U.K.}
\email{severini@cs.bris.ac.uk}
\date{June, 2002}
\thanks{Thanks to Sonia Mansilla for bring \cite{MS01} to my attention. Thanks to
Josef Lauri for our conversations on graph theory. }

\begin{abstract}
Study the relationship between unitary matrices and their patterns\ is
motivated by works in quantum chaology (see, \emph{e.g.}, \cite{KS99}) and
quantum computation (see, \emph{e.g.}, \cite{M96} and \cite{AKV01}). We
prove that if a Cayley digraph is a line digraph then it is the pattern of a
unitary matrix. We prove that for any finite group with two generators there
exists a set of generators such that the Cayley digraph with respect to such
a set is a line digraph and hence the pattern of a unitary matrix.
\end{abstract}

\maketitle

\section{Definitions}

The following standard definitions will be central:

\begin{definition}
Let $G=\left\langle S\right\rangle $ be a finite group. The (\emph{right}) 
\emph{Cayley digraph} $Cay\left( G,S\right) $ of $G$ with respect to $S$ is
the digraph whose vertex set is $G$, and whose arc set is the set of all
ordered pairs \{$\left( g,gs\right) :g\in G$ and $s\in S$\}.
\end{definition}

\begin{definition}
Let $D=\left( V,A\right) $ be a digraph. The \emph{line digraph} $L\left(
D\right) $ of $D$ is the digraph whose vertex set is $A\left( D\right) $,
and $\left( \left( u,v\right) ,\left( w,z\right) \right) \in A\left( L\left(
D\right) \right) $ if and only if $v=w$, where $u,v,w,z\in V\left( D\right) $
and $\left( u,v\right) ,\left( w,z\right) \in A\left( D\right) $.
\end{definition}

\begin{definition}
Let $M$ be a square matrix $M$ of size $n$. A digraph $D$ is said to be the 
\emph{pattern} of $M$, if $D$ is on $n$ vertices and, for every $u,v\in
V\left( D\right) $, $\left( u,v\right) \in A\left( D\right) $ if and only if
the entry $M_{uv}$ is nonzero.
\end{definition}

\begin{definition}
A square matrix $U$ with complex entries is said to be \emph{unitary} if it
is nonsingular and $U^{\dagger }U=I$, where $U^{\dagger }$ and $I$ denote
the adjoint of $U$ and the identity matrix, respectively.
\end{definition}

\section{On the Cayley digraphs that are patterns of unitary matrices}

Denote by $M\left( D\right) $ the adjacency matrix of a digraph $D$.

\begin{lemma}
\label{regular}If a regular digraph is a line digraph then it is the pattern
of a unitary matrix.
\end{lemma}

\begin{proof}
Recall the Richards characterization of line digraphs \cite{R67}:\ a digraph
is a line digraph if and only if: (i) the columns of its adjacency matrix
are identical or orthogonal; (ii) the rows of its adjacency matrix are
identical or orthogonal. Let $D$ be a regular line digraph. It is easy to
see that, from the Richards characterization, $M\left( D\right) $ is
composed of independent submatrices without zero entries. Since for any $n$
there exists a unitary matrix of size $n$ without zero entries, and since a
matrix composed by independent unitary matrices is unitary, the lemma is
proven.
\end{proof}

\begin{theorem}
\label{cayleyline}If a Cayley digraph is a line digraph then it is the
pattern of a unitary matrix.
\end{theorem}

\begin{proof}
By Lemma \ref{regular} and since a Cayley digraph is regular.
\end{proof}

\begin{remark}
The converse of Theorem \ref{cayleyline} is not true. Denote by $\mathbb{Z}%
_{n}$ the group of integers modulo $n$. The adjacency matrix of $\mathbb{Z}%
_{2}\times \mathbb{Z}_{2}=\left\langle \left( 0,0\right) ,\left( 1,0\right)
,\left( 0,1\right) \right\rangle $ is $\left[ 
\begin{array}{cccc}
1 & 1 & 1 & 0 \\ 
1 & 1 & 0 & 1 \\ 
1 & 0 & 1 & 1 \\ 
0 & 1 & 1 & 1
\end{array}
\right] $. Then, since the matrix $\frac{1}{^{\sqrt{3}}}\left[ 
\begin{array}{cccc}
1 & 1 & 1 & 0 \\ 
1 & -1 & 0 & 1 \\ 
-1 & 0 & 1 & 1 \\ 
0 & 1 & -1 & 1
\end{array}
\right] $ is real-orthogonal and hence unitary, $\mathbb{Z}_{2}\times 
\mathbb{Z}_{2}$ is the pattern of a unitary matrix, but it is not a line
digraph because it does not satisfy the Richard characterization.
\end{remark}

\begin{lemma}[\protect\cite{M01}, Proposition 4.3.1]
\label{sonia}If, for some $x\in S^{-1}$, $xS=H$, where $H$ is a subgroup of $%
G$ such that $\left| H\right| =r=\left| S\right| $, then $Cay\left(
G,S\right) =L\left( D\right) $, where $D$ is an $r$-regular (multi)digraph.
(See also \cite{MS01}, Corollary 7.)
\end{lemma}

Denote by $C_{n}$ the cyclic group of order $n$.

\begin{theorem}
\label{twogenerators}For any finite group with two generators, there exists
a set of generators such that the Cayley digraph with respect to such a set
is a line digraph and hence the pattern of a unitary matrix.
\end{theorem}

\begin{proof}
Let $G=\left\langle S\right\rangle $, where $S=\left\{ s_{1},s_{2}\right\} $%
. Take $s_{1}^{-1}\in S^{-1}$ (or, equivalently, $s_{2}^{-1}$). Then 
\begin{equation*}
s_{1}^{-1}S=\left\{ s_{1}^{-1}s_{1},s_{1}^{-1}s_{2}\right\} =\left\{
e,s_{1}^{-1}s_{2}\right\} .
\end{equation*}
Let $n$ be the order of $s_{1}^{-1}s_{2}$. Consider the group $%
C_{n}=\left\langle s_{1}^{-1}s_{2}\right\rangle $. Write $T=s_{1}C_{n}$.
Then $C_{n}=s_{1}^{-1}T$. Since $S\subset T$, $G=\left\langle T\right\rangle 
$, and since $\left| T\right| =n=\left| C_{n}\right| $ (in fact $T$ is a
left coset of $C_{n}$), by Lemma \ref{sonia}, $Cay\left( G,T\right) $ is a
line digraph, and hence, by Theorem \ref{cayleyline}, the pattern of a
unitary matrix.
\end{proof}

\begin{corollary}
For any finite simple group there exists a set of generators such that the
Cayley digraph of the group with respect to such a set is a line digraph and
hence the pattern of a unitary matrix.
\end{corollary}

\begin{proof}
By Theorem \ref{twogenerators}, together with the fact that every finite
simple group is generated by two elements \cite{AG84}.
\end{proof}

\section{Examples}

\begin{example}
By Theorem \ref{cayleyline}, $Cay\left( C_{n},g\right) $, for any $g\in
C_{n} $, is the only $1$-regular Cayley digraph that is the pattern of a
unitary matrix.
\end{example}

\begin{remark}
\label{remcay}By Theorem \ref{cayleyline}, a Cayley digraph $Cay\left(
G,\left\{ s_{1},s_{2}\right\} \right) $ is the pattern of a unitary matrix
if and only if $s_{1}=s_{2}s_{1}^{-1}s_{2}$ and $s_{2}=s_{1}s_{2}^{-1}s_{1}$%
. Moreover, if $G$ is abelian then $s_{1}^{2}=s_{2}^{2}$.
\end{remark}

Denote by $D_{n}$ the dihedral group of order $2n$ and by $e$ the identity
element.

\begin{example}
$Cay\left( D_{n},\left\{ s_{1},s_{2}\right\} \right) $ is the pattern of a
unitary matrix, because of Theorem \ref{cayleyline} and since the standard
presentation of $D_{n}$ is $\left\langle
s_{1},s_{2}:s_{1}^{n}=s_{2}^{2}=e,s_{2}s_{1}s_{2}=s_{1}^{-1}\right\rangle $
(see, e.g. \cite{CM72}, \S 1.5). Alternatively, by Theorem \ref{cayleyline}
and since 
\begin{equation*}
Cay\left( D_{n},\left\{ s_{1},s_{2}\right\} \right) \cong L\left( Cay\left( 
\mathbb{Z}_{n},\left\{ 1,-1\right\} \right) \right)
\end{equation*}
(see, e.g., \cite{BEFS95}).
\end{example}

\begin{example}
By Theorem \ref{cayleyline} and by Remark \ref{remcay}, $Cay\left( \mathbb{Z}%
_{4},\left\{ 1,-1\right\} \right) $ is the pattern of a unitary matrix.
\end{example}

\begin{definition}[\protect\cite{F84}]
Given the integers $k$ and $n$, $1\leq k\leq n-1$, denote by $P\left(
n,k\right) $ the digraph whose vertices are the permutations on $k$-tuples
from $\left[ n\right] =\left\{ 1,2,...,n\right\} $ and whose arcs are of the
form $\left( i_{1}\text{ }i_{2}\text{ ... }i_{k}\right) ,\left( i_{2}\text{ }%
i_{3}\text{ ... }i_{k}i\right) $, where $i\neq i_{1},i_{2},...,i_{k}$ (see
also \cite{BFF97}).
\end{definition}

Denote by $S_{n}$ the symmetric group on $\left[ n\right] =\left\{
1,2,...,n\right\} $.

\begin{example}
$Cay\left( S_{n},\left\{ s_{1},s_{2}\right\} \right) $, where $s_{1}=\left( 
\text{1 2 ... n}\right) $ and $s_{2}=\left( \text{1 2 ... n-1}\right) $, is
the pattern of a unitary matrix, by Theorem \ref{cayleyline} and since $%
Cay\left( S_{n},\left\{ s_{1},s_{2}\right\} \right) \cong L\left( P\left(
n,n-2\right) \right) $ (\cite{BFF97}, Lemma 2.1).
\end{example}

\begin{example}
$Cay\left( S_{n},T\right) $, where $T=\left( \text{1 2}\right) C_{n-1}$, is
the pattern of a unitary matrix. To show this, we use the construction in
the proof of Theorem \ref{twogenerators}. Consider $S_{n}=\left\langle
S=\left\{ \left( \text{1 2}\right) ,\left( \text{1 2 ... n}\right) \right\}
\right\rangle $ Then $S^{-1}=\left\{ \left( \text{1 2}\right) ,\left( \text{%
1 n ... 2}\right) \right\} $. Write $x=\left( \text{1 2}\right) \in S^{-1}$.
Then $\left( \text{1 2}\right) S=\left\{ e,\left( \text{2 3 ... n}\right)
\right\} $. Consider $C_{n-1}=\left\langle e,\left( \text{2 3 ... n}\right)
\right\rangle $. Write 
\begin{equation*}
T=\left( \text{1 2}\right) C_{n-1}=\left\{ \left( \text{1 2}\right) ,\left( 
\text{1 2}\right) \left( \text{2 3 ... n}\right) =\left( \text{1 2 ... n}%
\right) ,...\right\} .
\end{equation*}
Since $S\subset T$, $S_{n}=\left\langle T\right\rangle $. Moreover $x\in
T^{-1}$, $C_{n-1}=xT$ and $\left| T\right| =n-1=C_{n-1}$. Then, by Theorem 
\ref{cayleyline} and Lemma \ref{sonia}, $Cay\left( S_{n},T\right) $ is the
pattern of a unitary matrix.
\end{example}

\end{document}